\documentclass[12pt,reqno]{amsart}
\usepackage[T2A]{fontenc}
\usepackage[cp1251]{inputenc}
\usepackage[russian,english]{babel}

\usepackage{indentfirst}
\usepackage{verbatim}
\usepackage{amsbib}
\usepackage{amsfonts,amssymb,mathrsfs,amscd}
\usepackage{amsmath,latexsym}
\usepackage{amsthm}
\usepackage{verbatim}
\usepackage{eucal}
\usepackage{graphicx}
\usepackage{enumerate}
\usepackage{xcolor}

\theoremstyle{plain}

\newtheorem*{maintheorem}{Основная теорема}
\newtheorem{lemma}{Лемма}
\newtheorem{propos}{Предложение}
\newtheorem{corollary}{Следствие}

\renewcommand{\leq}{\leqslant} 
\renewcommand{\geq}{\geqslant}
\newcommand{\rad}{\text{\tiny\rm rad}}
\newcommand{\RR}{\mathbb{R}} 
\newcommand{\CC}{\mathbb{C}} 
\newcommand{\NN}{\mathbb{N}} 
\DeclareMathOperator{\Hol}{Hol} 
\DeclareMathOperator{\Exp}{Exp} 
\DeclareMathOperator{\Zero}{Zero} 
\DeclareMathOperator{\type}{type} 
 \DeclareMathOperator{\wid}{wid}
\DeclareMathOperator{\rh}{\text{\rm \tiny rh}}
 \DeclareMathOperator{\lh}{\text{\rm \tiny lh}}
\DeclareMathOperator{\dd}{\,{\mathrm d\!}}
\renewcommand{\Re}{{\rm Re \,}}
\renewcommand{\Im}{{\rm Im \,}}

\begin{document}

\title[Распределение нулей  целых функций экспоненциального типа \dots]{Распределение нулей  целых функций\\ экспоненциального типа  с ограничениями\\ на рост вдоль мнимой оси}

\author{А.\,Е.~Егорова, Б.\,Н.~Хабибуллин}

\selectlanguage{russian}
\maketitle

\section{Основные результаты}\label{s10}

Целая функция $f$ на \textit{комплексной плоскости\/} $\CC$ \textit{обращается в нуль\/} на последовательности точек ${\sf Z}=\{{\sf z}_k\}_{k=1,2,\dots}\subset \CC$, если для каждой точки $z\in \CC$ кратность нуля функции $f$ в этой точке  $z$ не меньше числа повторений этой точки $z$ в последовательности $\sf Z$ (пишем $f({\sf Z})=0$).Целую функцию $f$ называют \textit{целой функции экспоненциального типа} (пишем {\it ц.ф.э.т}), 
если конечен её тип
 \begin{equation}\label{typef}
\type_f:=\limsup_{z\to \infty}\frac{\ln |f(z)|}{|z|}.
\end{equation} 
 Если $f$ ---  ц.ф.э.т., то $2\pi$-периодическую функцию 
 \begin{equation*}
{\mathsf h}_f(\theta) :=\limsup_{r\to +\infty}\frac{\ln \bigl|f(re^{i\theta})\bigr|}{r}, \quad \theta\in \RR, 
\end{equation*}
называют её \textit{индикатором роста\/} \cite{Levin56}, \cite{Levin96}.

Последовательность ${\sf Z} =\{ {\sf z}_k\} \subset \CC$ {\it отделена от мнимой оси\/}
$i\RR$,  если 
\begin{equation}\label{con:dis}
\left(\liminf_{k\to\infty}\frac{\bigl| \Re {\sf z}_k \bigr|}{ |{\sf z}_k |} >0 \right)
\underset{\text{или}}{\Longleftrightarrow}
\left(\limsup_{k\to\infty}\frac{\bigl| \Im {\sf z}_k \bigr|}{ |{\sf z}_k |} <1 \right).
\end{equation}
Пара эквивалентных ограничений \eqref{con:dis} геометрически означает, что  найдётся  пара непустых открытых вертикальных углов, содержащих $i\RR\setminus \{0\}$, для которой точки ${\sf z}_k$ лежат вне этой пары углов  при всех достаточно больших $k$. 

Один из  результатов заметки --- это вытекающее из нашей основной теоремы
\begin{corollary}\label{thd} Пусть последовательность $\sf Z$ без предельных точек в $\CC$ отделена от мнимой оси, $b\in \RR^+:=\{x\in \RR \colon x\geq 0\}$ --- положительная полуось. 
Тогда эквивалентны  следующие два утверждения: 
\begin{enumerate}[{\rm I.}]
\item\label{fb1} Существует  ц.ф.э.т. $f\neq 0$, для которой $f({\sf Z})=0$, с индикатором роста
 \begin{equation}\label{indf}
{\mathsf h}_f(\pm \pi /2) \leq \pi b.
\end{equation}
\item\label{fb2} Существуют число $C\in \RR^+$ и  функция 
\begin{equation}\label{d}
d\colon \RR^+\to \RR^+, \quad \sup_{x\in \RR^+}d(x)<+\infty, \quad
\lim_{x\to +\infty}d(x)=0, 
\end{equation}
для которых 
\begin{equation}\label{Zb}
l_{\sf Z}(r,R):=
\max \left\{\sum_{\substack{r < |{\sf z}_k|\leq R\\\Re {\sf z}_k >0}} \Re \frac{1}{{\sf z}_k}, \sum_{\substack{r < |{\sf z}_k|\leq R\\\Re {\sf z}_k <0}} \Re \Bigl(-\frac{1}{{\sf z}_k}\Bigr)\right\}
\leq \bigl(b+d(R)\bigr) \ln \frac{R}{r} +C 
\end{equation}
при всех $0< r<R<+\infty$.
\end{enumerate}
\end{corollary}
В частном случае для \textit{положительной последовательности\/}   ${\sf Z}$, расположенной исключительно на \textit{положительной полуоси\/} 
$\RR^+$,  наше следствие  \ref{thd} --- это в точности  один из классических совместных результатов  П.~Мальявена и Л.~А.~Рубела \cite[теорема 6.3]{MR}. В \cite[следствие 4.1]{kh91AA},  \cite[теорема 2.5.2]{KhDD92} и \cite[теорема 3.2.5, следствие 3.2.1]{Khsur} следствие \ref{thd} было 
анонсировано без доказательства. Мы выводим его  из нашей основной теоремы  с полным доказательством. Несколько неожиданно для нас это полное доказательство оказалось не столь уж и простым, во всяком случае, в нашей реализации, несмотря на то что первоначально при его анонсировании оно представлялось несложным переносом близкой схемы П.~Мальявена   и Л. А. Рубела из \cite[доказательства теорем 6.3--5]{MR} в реалии основного результата из \cite[основная теорема]{kh91AA}.
Наша основная теорема ниже приводится в форме, продиктованной формулировкой основного результата из  \cite[теорема 4.1]{MR}, которому  посвящен один из основных разделов  монографии Л.~А.~Рубела  с Дж.~Э.~Коллиандром \cite[22, основная теорема]{RC}, а также формулировками близких, но иных результатов из  \cite[основная теорема]{KhaD88}, \cite[основная теорема]{Kha89}, \cite[основная теорема]{kh91AA}, 
\cite[теоремы 2.4.1, 2.4.2]{KhDD92}, \cite[теорема 3.2.1]{Khsur}.

Для $S\subset \CC$ через $\Hol(S)$ обозначаем \textit{векторное пространство над $\CC$ всех голоморфных функций\/} в какой-либо своей открытой окрестности множества $S$.

Следуя \cite[2]{MR} и  \cite[22]{RC}, последовательности  ${\sf Z}=\{{\sf z}_k\}\subset \CC$ 
сопоставляем идеал 
\begin{equation*}
I({\sf Z}):=\bigl\{f\in \Hol (\CC)\colon f({\sf Z})=0\bigr\}\subset \Hol (\CC)
\end{equation*} 
в кольце $\Hol (\CC)$, а также идеал в кольце всех ц.ф.э.т.\footnote{В \cite[2]{MR} и  \cite[22]{RC} идеал $I^1({\sf Z})$ обозначен соотв. как $\mathcal F ({\sf Z})$ и $F({\sf Z})$.}
\begin{equation*}
I^1({\sf Z}):=I({\sf Z})\cap \bigl\{f\in \Hol (\CC)\colon {\type}_f\overset{\eqref{typef}}{<}+\infty \bigr\}.
\end{equation*} 
Полагаем $\Hol_*(\CC):=\Hol(\CC)\setminus \{0\}$ и 
\begin{equation*}
I_*({\sf Z}):=I({\sf Z})\setminus \{0\}=I({\sf Z})\cap \Hol_*(\CC),
\quad 
I_*^1({\sf Z}):=I^1({\sf Z})\setminus \{0\}=I^1({\sf Z})\cap \Hol_*(\CC).
\end{equation*}

Последовательности ${\sf Z}\subset \CC$ без предельных точек в $\CC$ сопоставляем \textit{считающие меру\/} и \textit{функцию\/} соответственно 
\begin{equation}\label{nZ}
n_{\sf Z}(S):=\sum_{{\sf z}_k\in S}1, \; S\subset \CC, 
 \qquad n_{\sf Z}^{\rad}(r):=\sum_{|{\sf z}_k|\leq r} 1, \; r\in \RR^+.
\end{equation}
Последовательность ${\sf Z}$  {\it конечной верхней плотности,\/}
если 
\begin{equation}\label{typeZ}
{\type}_{\sf Z}:=\limsup_{0<r\to +\infty}\frac{n_{\sf Z}^{\rad}(r)}{r}<+\infty.
\end{equation}
Последовательность конечной верхней плотности не имеет предельных точек в $\CC$.
\begin{propos}[{\rm \cite{Levin56},\cite{Levin96}, \cite{RC}}]\label{pr2_1}
 $I_*({\sf Z})\neq \varnothing$, если и только если   $\sf Z$ не имеет предельных точек в\/ $\CC$ --- {\rm частный случай классической теоремы Вейерштрасса для $\CC$.}  $I^1_*({\sf Z})\neq \varnothing$, если и только если $\sf Z$ конечной верхней плотности  ---{\rm  частный случай классической теоремы Адамара\,--\,Вейерштрасса о представлении целых функций.} 
\end{propos}

Последовательность всех нулей целой функции $g\in \Hol_*(\CC)$, перенумерованную каким-либо образом с учётом кратности, обозначаем через  $\Zero_g$.

Через $\CC_{\rh}:=\{z\in \CC\colon \Re z>0\}$ обозначаем \textit{открытую правую полуплоскость}.

\begin{maintheorem} Пусть две последовательности ${\sf Z}\subset \CC$ и ${\sf W}\subset \CC_{\rh}$ конечной верхней плотности  отделены от $i\RR$. 
Тогда эквивалентны три утверждения: 
\begin{enumerate}[{\rm I.}]
\item\label{fgi} Для любой функции $g\in I_*^1({\sf W})$ с последовательностью нулей, отделённой от $i\RR$,  найдётся функция $f\in I^1_*({\sf Z})$ с ограничением 
\begin{equation}\label{fgiR}
\ln \bigl|f(iy)\bigr|\leq \ln \bigl|g(iy)\bigr|+o(|y|)\quad \text{при  $y\to  +\infty$}.
\end{equation}
\item\label{fgiii} Существует пара функций $f\in I_*^1({\sf Z})$ и  $g\in I^1_*({\sf W})$ c $\Zero_g \cap\, \CC_{\rh} ={\sf W}$, т.\,е. с сужением  $n_{\Zero_g}\bigm|_{\CC_{\rh}}=n_{\sf W}$,  для которой имеет место   \eqref{fgiR}. 
\item\label{fgii} Существуют $C\in \RR^+$ и   функция $d$ из \eqref{d}, для которых в обозначении \eqref{Zb}
\begin{equation}\label{Zld}
l_{\sf Z}(r,R)\leq l_{\sf W}(r,R) +d(R)\ln \frac{R}{r}
+C\quad \text{при всех\/ $0< r<R<+\infty$}.
\end{equation}
\end{enumerate}

\end{maintheorem}

Для функций $v\colon i\RR \to \RR_{\pm \infty}:=\{-\infty\}\cup \RR\cup \{+\infty\}$ определим  интеграл 
\begin{equation}\label{fK:abp+}
J_{i\RR}(r,R;v):=
\frac{1}{2\pi}\int_r^{R} \frac{v(-iy)+v(iy)}{y^2} \dd y, \quad 0<r<R<+\infty.
\end{equation}

Следующий результат более общий, чем следствие  \ref{thd}.

\begin{corollary}\label{cor2} Пусть как последовательность 
нулей  ц.ф.э.т. $g\neq 0$, так и последовательность ${\sf Z}$ 
отделены от $i\RR$. Эквивалентны два  утверждения:    
\begin{enumerate}[{\rm I.}]
\item\label{gb1} Существует  ц.ф.э.т. $f\in I_*^1({\sf Z})$, для которой выполнено 
\eqref{fgiR}. 
\item\label{gb2} Существуют число 
$C\in \RR^+$ и  функция $d$ из \eqref{d}, для которых 
\begin{equation}\label{gb}
l_{\sf Z}(r,R)\leq J_{i\RR}\bigl(r,R;\ln|g|\bigr)+d(R) \ln \frac{R}{r} +C 
\quad \text{при всех\/ $0< r<R<+\infty$}.
\end{equation}
\end{enumerate}
\end{corollary}
Следствие  \ref{thd}  сразу  получается из следствия \ref{cor2}
при выборе $g(z):=\sin \pi b z$, $z\in \CC$.

Следующий результат имеет форму теоремы о мультипликаторе.
\begin{corollary}\label{cormult} 
Пусть $p, g$ --- две ц.ф.э.т.  с последовательностями нулей, отделёнными от мнимой оси. Следующие два утверждения эквивалентны:
\begin{enumerate}[{\rm I.}]
\item\label{hb1} Существует  ц.ф.э.т.-мультипликатор 
 $h\neq 0$, для которой ц.ф.э.т.-про\-и\-з\-в\-е\-д\-е\-н\-ие $f:=ph$ удовлетворяет условию  \eqref{fgiR}. 
\item\label{pb2} Существуют число $C\in \RR^+$ и  функция $d$ из \eqref{d},  для которых 
\begin{equation}\label{pb}
J_{i\RR}\bigl(r,R;\ln |p|\bigr)\leq J_{i\RR}\bigl(r,R;\ln|g|\bigr)+d(R) \ln \frac{R}{r} +C 
\end{equation}
при всех\/ $0< r<R<+\infty$.
\end{enumerate}
\end{corollary}

\textit{Шириной множества $S\subset \CC$ в направлении} $0$
называется число 
\begin{equation*}
{\wid}_0 (S):=\sup\{\Im z_1-\Im z_2 \colon z_1,z_2\in S \}, 
\end{equation*}
которое в  терминах $2\pi$-периодической опорной функции 
\begin{equation}\label{sS}
{\sf s}_S(\theta):=\sup_{z\in S} \Re ze^{-i\theta}, \quad \theta\in \RR,
\end{equation}
 выражается как 
\begin{equation}\label{ws}
{\wid}_0 (S)={\sf s}_S(-\pi/2)+{\sf s}_S(\pi/2).
\end{equation}

Через $\Hol'(\CC)$ обозначаем векторное пространство над $\CC$
\textit{аналитических функционалов\/} \cite[4.7]{Ho}, \cite[\S~4, 5]{Nap82}, т.\,е. всех линейных непрерывных функционалов на пространстве $\Hol(\CC)$, снабжённом  топологией равномерной сходимости на компактах.  Для аналитического функционала 
$\mu\in \Hol'(\CC)$  его \textit{преобразование Лапласа\/}
\begin{equation}\label{FL}
\widehat \mu (z):=\int_{\CC} e^{zw} \dd \mu(w), \quad z\in \CC,
\end{equation}
 --- целая функция экспоненциального типа. Для пары аналитических функционалов  $\mu,\nu\in \Hol'(\CC)$ определена  операция \textit{свёртки\/} $\mu\ast\nu\in \Hol'(\CC)$,  действующая по правилу
\begin{equation*}
(\mu\ast \nu)(f)=\mu_w\bigl(\nu_z (f(z+w))\bigr), \quad f\in \Hol(\CC).
\end{equation*}
Компакт $K\subset \CC$  \textit{определяющий для\/} $\mu\in \Hol'(\CC)$, если 
для каждой открытой окрестности $O$ компакта $K\subset O$ найдётся постоянная 
$C_O\in \RR^+$, с которой 
\begin{equation*}
\bigl|\mu(f)\bigr|\leq C_O\sup_O|f|\quad\text{для всех $f\in \Hol(\CC)$}.
\end{equation*} 

\begin{corollary}\label{coraf} Пусть $0\neq \mu\in \Hol'(\CC)$ и последовательность нулей 
ц.ф.э.т. $\widehat \mu$ отделена от нуля, $b\in \RR^+$. 
 Эквивалентны два утверждения:
\begin{enumerate}[{\rm I.}]
\item\label{mub1} Существует  ненулевой  функционал 
 $\nu\in \Hol'(\CC)$, для которого свёртка $\mu\ast \nu$ обладает определяющим (выпуклым) компактом $K$ ширины ${\wid}_0(K)\leq 2\pi b$.
\item\label{mub2} Существуют число $C\in \RR^+$ и  функция $d$ из \eqref{d},  для которых 
\begin{equation}\label{mub}
J_{i\RR}\bigl(r,R;\ln |\widehat \mu|\bigr)\leq \bigl(b+d(R) \bigr)\ln \frac{R}{r} +C 
\quad\text{при всех\/ $0< r<R<+\infty$.}
\end{equation}
\end{enumerate}
\end{corollary}

Через $\NN :=\{1,2, \dots\}$ обозначаем \textit{множество натуральных чисел,}
$\NN_0:=\{0\}\cup \NN$. 

Последовательности ${\sf Z}=\{{\sf z}_k\}_{k=1,2,\dots}\subset \CC$ сопоставляем 
{\it экспоненциальную систему $\Exp^{\sf Z}\subset \Hol_* (\CC)$
с последовательностью показателей\/} $\sf Z$: 
\begin{equation}\label{ExpZ} 
\Exp^{\sf Z}:=\bigl\{ z\mapsto z^pe^{{\sf z}_kz}
\colon z\in \CC,\;  p\in \NN_0, \;0\leq p\leq n_{\sf Z}\bigl(\{{\sf z}_k\}\bigr)-1\bigr\}.
\end{equation} 
Пространство $\Hol(K)$ локально аналитических функций на компакте $K\subset\CC$ наделяется  \textit{топологией индуктивного предела\/} \cite{Seb57}, \cite[гл.~XI, \S~5]{KA},  \cite[\S~4,~4]{Nap82}, \cite[0.1.4]{Khsur}.

\begin{corollary}\label{corcompl} Пусть последовательность  ${\sf Z}$ 
отделена от $i\RR$, $b\in \RR^+$. Эквивалентны два утверждения:
\begin{enumerate}[{\rm I.}]
\item\label{Eub1} Система $\Exp^{\sf Z}$ полна в каждом  пространстве 
$\Hol(K)$, когда $K$ --- выпуклый компакт ширины ${\wid}_0(K)\leq 2\pi b$. 
\item\label{Eub2} Справедливо соотношение 
\begin{equation}\label{Eub}
\inf_{d}\sup_{1\leq r<R<+\infty}
\Bigl(l_{\mathsf Z}(r,R)-\bigl(b+d(R)\bigr)\ln \frac{R}{r}\Bigr)=+\infty,
\end{equation}
где $\inf$ в левой части берётся  по  всех функциям  $d$ вида  \eqref{d}.
\end{enumerate}
\end{corollary}

\section{Доказательство основной теоремы}\label{mth}
\setcounter{equation}{0}
Выбор ц.ф.э.т.
\begin{equation*}
g(z):=\prod_{k} \Bigl(1-\frac{z^2}{{\sf w}_k^2}\Bigr), 
\quad {\mathsf W}=\{{\sf w}_k\}_{=1,2,\dots}\subset \CC_{\rh}, 
\quad z\in \CC,
\end{equation*} 
в п. \ref{fgiii} при выполнении \ref{fgi} доказывает импликацию \ref{fgi}$\Rightarrow$\ref{fgiii}.

\begin{proof}[{\bf Доказательство импликации \ref{fgiii}$\Rightarrow$\ref{fgii}}] Интегрируя  
как в \eqref{fK:abp+} условие  \eqref{fgiR},  для некоторого фиксированного 
числа $r_0>0$  получаем 
\begin{equation*}
J_{i\RR}\bigl(r,R;\ln |f|\bigr)\leq J_{i\RR}\bigl(r,R;\ln |g|\bigr)+
\int_r^R\frac{Q(y)}{y^2}\dd y+C\quad\text{для всех $r_0\leq r<R<+\infty$},
\end{equation*}   
где ограниченная функция удовлетворяет условию \eqref{diy}. 

\begin{lemma}[{\cite[Следствие]{Kha19}}]\label{corln} Пусть $0<r_0\in \RR^+$. 
Если  для положительной функции $Q\colon [r_0,+\infty) \to \RR^+$ существует предел 
\begin{equation}\label{diy}
\lim_{x\to +\infty}\frac{Q(x)}{x}=0,
\end{equation}
то найдётся  убывающая функция   $d\colon [r_0,+\infty)\to \RR^+$, для которой 
\begin{subequations}\label{dQ}
\begin{align}
\int_r^R\frac{Q(x)}{x^2}\dd t&\leq d(R)\ln\frac{R}{r}\quad \text{при всех $r_0\leq r<R<+\infty$},
\tag{\ref{dQ}A}\label{{ad}A}\\
\lim_{R\to +\infty} d(R)&=0. 
\tag{\ref{dQ}$_0$}\label{{ad}0}
\end{align}
\end{subequations} 
Если для функции $d\colon [r_0,+\infty) \to \RR^+$ выполнено \eqref{{ad}0}, то найдётся возрастающая функция $Q\colon [r_0,+\infty)\to \RR^+$, для которой выполнено \eqref{diy} и 
\begin{equation}\label{Qd}
d(R)\ln\frac{R}{r} \leq \int_r^R\frac{Q(x)}{x^2}\dd x
\quad \text{при всех $r_0\leq r<R<+\infty$}.
\end{equation}
\end{lemma}

По первой части леммы \ref{corln} это неравенство можно переписать как 
\begin{equation}\label{Jl}
J_{i\RR}\bigl(r,R;\ln |f|\bigr)\leq J_{i\RR}\bigl(r,R;\ln |g|\bigr)+
d(R)\ln \frac{R}{r}+C\quad\text{для всех $r_0\leq r<R<+\infty$},
\end{equation}   
где $d$ --- некоторая функция вида \eqref{d}. Положим 
\begin{equation}\label{Zbrl}
l_{\sf Z}^{\rh}(r,R):=\sum_{\substack{r < |{\sf z}_k|\leq R\\\Re {\sf z}_k >0}} \Re \frac{1}{{\sf z}_k}, \quad
l_{\sf Z}^{\lh}(r,R):=\sum_{\substack{r < |{\sf z}_k|\leq R\\\Re {\sf z}_k <0}} \Re \Bigl(-\frac{1}{{\sf z}_k}\Bigr),
\end{equation}
что позволяет записывать  левую часть \eqref{Zb} как
\begin{equation}\label{lllr}
l_{\sf Z}(r,R)=\max \bigl\{l_{\sf Z}^{\lh}(r,R), l_{\sf Z}^{\rh}(r,R)\bigr\}.
\end{equation}

\begin{lemma}[{\cite[(1.3)]{Kha89}, \cite[(0.4)]{kh91AA}, \cite[предложение 4.1, (4.19)]{KhII}}]\label{lemJl} При фиксированном числе $r_0>0$ для любой ц.ф.э.т. $f\neq 0$ имеет место соотношение
\begin{equation*}
\sup_{r_0\leq r<R<+\infty}\max \Bigl\{\bigl|J_{i\RR}(r,R;\ln |f|)-l_{\Zero_f}^{\rh}(r,R)\bigr|,\;
\bigl|J_{i\RR}(r,R;\ln|f|)-l_{\Zero_f}^{\lh}(r,R)\bigr|\Bigr\}
<+\infty.
\end{equation*}
\end{lemma}
 Из \eqref{Jl} по лемме \ref{lemJl}, применённой и к $f$, и к $g$,  в обозначении \eqref{lllr} получаем
\begin{multline}\label{Zdg}
l_{\sf Z}(r,R)\leq l_{\Zero_f}(r,R)\leq
 l_{\Zero_g}^{\rh}(r,R) +d(R)\ln \frac{R}{r}+C'\\
= l_{\sf W}(r,R) +d(R)\ln \frac{R}{r}
+C'\quad \text{при всех\/ $r_0\leq  r<R<+\infty$}.
\end{multline}
с постоянной $C'$, не зависящей от $r$ и $R$. По определениям \eqref{Zbrl}
при достаточно малом  $r_0>0$ его в \eqref{Zdg} можно заменить на $0$, но с $r>0$, что даёт требуемое \eqref{Zld}. 
\end{proof}

\begin{proof}[{\bf Доказательство импликации \ref{fgii}$\Rightarrow$\ref{fgi}}] По второй части леммы \ref{corln} неравенства  \eqref{Zld} можно переписать для сколь угодно малого 
$r_0>0$ как 
\begin{equation}\label{Zld+}
l_{\sf Z}(r,R)\leq l_{\sf W}(r,R) +
\int_r^R\frac{Q(y)}{y^2}\dd y+C\quad\text{для всех $r_0\leq r<R<+\infty$},
\end{equation}   
где $Q\colon \RR^+\to \RR^+$ --- \textit{возрастающая функция,\/} удовлетворяющая предельному соотношению \eqref{diy}. Рассмотрим  последовательность 
${\mathsf Q}:=\{{\sf q}_k\}_{k=1,2,\dots}\subset \RR^+$, однозначно определяемую через её  считающую функцию 
\begin{equation}\label{nQ}
n_{\sf Q}^{\rad}(r)\overset{\eqref{nZ}}{:=}\lfloor Q(r)\rfloor:=
\max\bigl\{m\in \NN_0\colon m\leq Q(r)\bigr\} \text{ --- \textit{целая часть} $Q(r)$}, \; r\in \RR^+. 
\end{equation}
Для такой последовательности ${\sf Q}$, очевидно,
\begin{equation}\label{typeq}
\type_{\sf Q}\overset{\eqref{typeZ}}{\leq} \limsup_{r\to +\infty}\frac{Q(r)}{r}\overset{\eqref{diy}}{=}0, 
\end{equation}
а также по определениям \eqref{Zbrl}--\eqref{lllr}
\begin{equation*}
l_{\sf Q}(r,R)=\int_r^R\frac{1}{t}\dd n_{\sf Q}^{\rad}(t)
\overset{\eqref{nQ}}{=}\frac{\lfloor Q(R)\rfloor}{R}-\frac{\lfloor Q(r)\rfloor}{r}+
\int_r^R\frac{\lfloor Q(R)\rfloor}{y^2}\dd y,
\end{equation*}
откуда 
\begin{equation*}
\int_r^R\frac{Q(y)}{y^2}\dd y \leq 
l_{\sf Q}(r,R)+\int_r^R\frac{\dd y}{y^2}+C'\leq l_{\sf Q}(r,R)+C_0, 
\end{equation*}
где числа $C', C_0$ не зависят от значений  $R>r\geq r_0$. Таким образом,  \eqref{Zld+} влечёт за собой существование числа $C_1\in \RR^+$, с которым
\begin{equation}\label{WQ}
l_{\sf Z}(r,R)\leq l_{{\sf W}\cup{\sf Q}}(r,R) +
C\quad\text{для всех $r_0\leq r<R<+\infty$},
\end{equation}  
где объединение ${\sf W}\cup{\sf Q}$ двух последовательностей --- эта последовательность со считающей мерой $n_{{\sf W}\cup{\sf Q}}:=
n_{\sf W}+n_{\sf Q}$. Тем более, после дополнения  последовательности ${\sf W}\cup{\sf Q}$ новой последовательностью $-{\sf Q}:=\{{-\sf q}_k\}_{k=1,2,\dots}\subset -\RR^+$ сохраняется \eqref{WQ}:
\begin{equation}\label{WQ-}
l_{\sf Z}(r,R)\leq l_{{\sf W}\cup{\sf Q}\cup ({-\sf Q})}(r,R) +
C\quad\text{для всех $r_0\leq r<R<+\infty$},
\end{equation}
Пусть $g\in I_*^1({\sf W})$ --- произвольная ц.ф.э.т., обращающаяся в нуль на  $\sf W$ с последовательностью нулей $\Zero_g\supset {\sf W}$, отделённой от $i\RR$. Из \eqref{WQ-}
следует
\begin{equation}\label{WQ-g}
l_{\sf Z}(r,R)\leq l_{{\Zero_g}\cup{\sf Q}\cup ({-\sf Q})}(r,R) +
C\quad\text{для всех $r_0\leq r<R<+\infty$}.
\end{equation}
Рассмотрим ц.ф.э.т. 
\begin{equation}\label{q}
q(z):=\prod_{k} \Bigl(1-\frac{z^2}{{\sf q}_k^2}\Bigr), 
\quad z\in \CC.
\end{equation} 
По построению $\Zero_q={\sf Q}\cup ({-\sf Q})$. Таким образом, из 
\eqref{WQ-g} получаем
\begin{equation*}
l_{\sf Z}(r,R)\leq l_{{\Zero_g}\cup{\Zero_q}}(r,R) +
C=l_{\Zero_{gq}}(r,R) +C
\quad\text{для всех $r_0\leq r<R<+\infty$},
\end{equation*}
где произведение  $gq$ --- ц.ф.э.т. с последовательностью нулей, отделённой от $i\RR$. Отсюда, применяя лемму \ref{lemJl} к ц.ф.э.т. 
 $gq$,  имеем 
\begin{equation*}
l_{\sf Z}(r,R)\leq J_{i\RR}\bigl(r,R; \ln|gq|\bigr) +C
\quad\text{для всех $r_0\leq r<R<+\infty$}.
\end{equation*}
Из последнего  по \cite[основная теорема]{kh91AA}  следует, что найдётся ц.ф.э.т. $f\neq 0$, обращающаяся в нуль на ${\sf Z}$, для которой выполнено неравенство 
\begin{equation*}
\ln \bigl|f(iy)\bigr|\leq \ln \bigl|(gq) (iy)\bigr|=
\ln \bigl|g (iy) \bigr|+\ln\bigl|q(iy)\bigr|\quad\text{для всех $y\in \RR$}.
\end{equation*}
Здесь для ц.ф.э.т. $q$ из \eqref{q} с симметричной относительно нуля 
последовательностью нулей с нулевой верхней плотностью \eqref{typeq} по теореме Адамара \cite[теорема 15]{Levin56} имеем $\ln\bigl|q(iy)\bigr|\leq o\bigl(|y|\bigr)$
при $|y|\to +\infty$, что даёт \eqref{fgiR}. 
\end{proof}

\section{Доказательства следствий из основной теоремы}
\setcounter{equation}{0}
\begin{proof}[{\bf Доказательство следствия \ref{cor2}}] 
Пусть ${\sf W}:=\Zero_g\cap \CC_{\rh}$. Если выполнено утверждение 
 \ref{gb1} следствия \ref{cor2}, то имеет место утверждение \ref{fgiii} основной теоремы. 
Следовательно,  по импликации \ref{fgiii}$\Rightarrow$\ref{fgii} основной теоремы 
имеет место утверждение \ref{fgii} основной теоремы, где в правой части соотношения  
\eqref{Zld}  величину-характеристику
\begin{equation}\label{lW}
l_{\sf W}(r,R)\overset{\eqref{Zbrl}}{=}l_{\sf W}^{\rh}(r,R)=l_{\Zero_g}^{\rh}(r,R)
\end{equation}
по лемме \ref{lemJl} можно заменить на $J_{i\RR} \bigl(r,R;\ln |g|\bigr)$, что доказывает 
импликацию  \ref{gb1}$\Rightarrow$\ref{gb2} следствия \ref{cor2}.
Обратно, если выполнено \ref{gb2}, то вновь  по лемме \ref{lemJl} в правой части \eqref{gb} можно заменить  $J_{i\RR} \bigl(r,R;\ln |g|\bigr)$ на характеристику \eqref{lW} и, в частности,  для некоторого числа $C\in \RR^+$
\begin{equation*}
l_{\sf Z}(r,R)\leq C\ln \frac{R}{r} +C\quad \text{для всех $0<r<R<+\infty$}.
\end{equation*} 
Отсюда ввиду отделённости от  мнимой оси  ${\sf Z}$ --- последовательность конечной верхней плотности.
Таким образом,  выполнено  утверждение \ref{fgii} основной теоремы и из импликации \ref{fgii}$\Rightarrow$\ref{fgi} основной теоремы получаем утверждение \ref{gb1} следствия \ref{cor2}.  
\end{proof}

\begin{proof}[{\bf Доказательство следствия \ref{cormult}}] Утверждение \ref{hb1} следствия 
\ref{cormult} эквивалентно утверждению 
 \ref{gb1} следствия \ref{cor2} для случая ${\sf Z}:=\Zero_p$, что по следствию 
\ref{cor2} эквивалентно утверждению \ref{gb2} следствия \ref{cor2}, где в неравенствах 
\eqref{gb} по лемме \ref{lemJl} характеристику $l_{\sf Z}(r,R)=l_{\Zero_p}(r,R)$ 
можно заменить на $J_{i\RR} \bigl(r,R;\ln |p|\bigr)$ как в  \eqref{pb}. 
\end{proof}

\begin{proof}[{\bf Доказательство следствия \ref{coraf}}] Утверждение \eqref{mub1}
 эквивалентно тому, что для ц.ф.э.т. $p:=\widehat{\mu}\neq 0$ существует ц.ф.э.т.-мультипликатор $h=\widehat \nu\neq 0$, для которой индикаторная (или сопряжённая) диаграмма  \cite{Levin56}, \cite{Levin96} ц.ф.э.т.-произведения $q:=ph$ ширины  не более $2\pi b$ в направлении $0$. В терминах индикатора  роста    ${\sf h}_q$ это эквивалентно неравенству ${\sf h}_q(\pi/2)+{\sf h}_q(-\pi/2)\leq 2\pi b$. Тогда можно подобрать такое вещественное число $a\in \RR$, что индикатор роста ${\sf h}_f$
ц.ф.э.т. $f(z):=e^{iaz}q(z)$, $z\in \CC$,  удовлетворяет неравенствам ${\sf h}_f(\pm\pi/2)\leq \pi b$. Это означает, что с функцией 
 $g_b(z)=\sin \pi bz$, $z\in \CC$,  выполнено неравенство 
\begin{equation*}
\ln |p(iy)|+\ln|h(iy)e^{-ay}|=\ln |q(iy)e^{-ay}|=\ln |f(iy)|\leq \ln |g_b(iy)|+o\bigl(|y|\bigr), 
\quad |y|\to +\infty.  
\end{equation*}
Таким образом, выполнено утверждение \ref{hb1} следствия \ref{cormult}. Тогда по эквивалентному ему утверждению \ref{pb2} следствия \ref{cormult} 
из  неравенств \eqref{pb} имеем эквивалентные им при  $g:=g_b$ неравенства 
\begin{multline*}\label{pb+}
J_{i\RR}\bigl(r,R;\ln |\widehat \mu|\bigr)=J_{i\RR}\bigl(r,R;\ln |p|\bigr)\leq J_{i\RR}\bigl(r,R;\ln|g_b|\bigr)+d(R) \ln \frac{R}{r} +C \\
\leq b\ln\frac{R}{r}+d(R) \ln \frac{R}{r} +C'
\quad\text{при всех\/ $0< r<R<+\infty$.}
\end{multline*}

\end{proof}

\begin{proof}[{\bf Доказательство следствия \ref{corcompl}}] Сильное сопряжённое к пространству $\Hol(K)$, когда $K$ --- выпуклый компакт  в $\CC$ с опорной функцией ${\sf s}_K$ из \eqref{sS}, состоит из всех аналитических функционалов $\mu\in \Hol'(\CC)$, для которых $K$ --- определяющее множество. Множество  преобразований Лапласа  $\widehat \mu$  \eqref{FL} всех таких аналитических функционалов $\mu$  --- это в точности пространство всех ц.ф.э.т. $f$ с индикатором роста ${\sf h}_f(\theta)\leq {\sf s}_K(-\theta)$, $\theta\in \RR$ \cite[теорема 4.7.11]{Ho}. По теореме Хана\,--\,Банаха система $\Exp^{\sf Z}$ не полна в пространстве $\Hol (K)$ тогда и только тогда, когда существует ненулевой аналитический функционал $\mu$ с определяющим множеством $K$, который принимает значение $0$ на всех функциях из    $\Exp^{\sf Z}$. В терминах их преобразований Лапласа   \eqref{FL} это означает, что существует ц.ф.э.т. $f_K\neq 0$ с $f_K({\sf Z})=0$ и с индикатором роста ${\sf h}_{f_K}(\theta)\leq {\sf s}_K(-\theta)$, $\theta\in \RR$. Таким образом, система $\Exp^{\sf Z}$ {\it не полна в  каком-нибудь\/}  пространстве  $\Hol(K)$, где $K$ --- {\it выпуклый компакт ширины  не больше $2\pi b$ в направлении  $0$,\/} если и только если  указанная ц.ф.э.т. $f_K$ удовлетворяет условию 
\begin{equation*}
{\sf h}_{f_K}(\pi/2)+{\sf h}_{f_K}(-\pi/2)\leq 
{\sf s}_K(-\pi/2)+{\sf s}_K(\pi/2)\overset{\eqref{ws}}{\leq} 2\pi b.
\end{equation*} 
Последнее эквивалентно существованию ц.ф.э.т. $f$  с $f({\sf Z})=0$ и с индикатором роста ${\sf h}_f(\pm\pi/2)\leq \pi b$, что по следствию \ref{thd} эквивалентно выполнению утверждения \ref{fb2} следствия \ref{thd}. Таким образом, система   
$\Exp^{\sf Z}$ {\it полна в  любом\/}  пространстве  $\Hol(K)$, где $K$ --- выпуклый компакт ширины  не больше $2\pi b$ в направлении  $0$, если и только если имеет место отрицание утверждения  \ref{fb2} следствия \ref{thd}.  Это отрицание и есть  утверждение 
\ref{Eub2} следствия   \ref{corcompl}.
\end{proof}

\end{document}